\documentclass[conference]{IEEEtran}
\IEEEoverridecommandlockouts
\usepackage{cite}
\usepackage{amsmath,amssymb,amsfonts}
\usepackage{graphicx}
\usepackage{textcomp}
\usepackage{xcolor}
\usepackage{graphics} 
\usepackage{epsfig}
\usepackage{mathptmx}
\usepackage{times} 
\usepackage{amsmath} 
\usepackage{amssymb}  
\usepackage[nocomma]{optidef}
\usepackage{xcolor}
\usepackage{bm}
\usepackage{booktabs} 
\usepackage{algorithm}  
\usepackage{algorithmicx} 
\usepackage{cite}
\usepackage{multirow}
\usepackage{threeparttable}
\usepackage{subfigure}
\usepackage{arydshln}

\def\BibTeX{{\rm B\kern-.05em{\sc i\kern-.025em b}\kern-.08em
    T\kern-.1667em\lower.7ex\hbox{E}\kern-.125emX}}

\begin{document}
\title{Mitigation-Aware Bidding Strategies in Electricity Markets
\thanks{JA is partially supported by NSF grants EECS-2144634 \& ECCS-2231350.
YW and JK are partially supported by the U.S. Department of Energy's Office of Energy Efficiency and Renewable Energy under the Solar Energy Technology Office Award Number DE-EE0008769. The views expressed herein do not necessarily represent the views of the U.S. Department of Energy or the U.S. Government.
}
}

\author{\IEEEauthorblockN{ Yiqian Wu}
\IEEEauthorblockA{\textit{Department of Electrical Engineering} \\
\textit{Columbia University}\\
New York, NY 10027, USA \\
yw3740@columbia.edu}
\and
\IEEEauthorblockN{Jip Kim}
\IEEEauthorblockA{\textit{Department of Energy Engineering} \\
\textit{Korea Institute of Energy Technology}\\
Naju, South Jeolla, Korea \\
jipkim@kentech.ac.kr}
\and
\IEEEauthorblockN{James Anderson}
\IEEEauthorblockA{\textit{Department of Electrical Engineering} \\
\textit{Columbia University}\\
New York, NY 10027, USA \\
james.anderson@columbia.edu}
}

\maketitle

\begin{abstract}
Market power exercise in the electricity markets distorts market prices and diminishes social welfare.
Many markets have implemented market power mitigation processes to eliminate the impact of such behavior.
The design of mitigation mechanisms has a direct influence on investors' profitability and thus mid-/long-term resource adequacy.
In order to evaluate the effectiveness of the existing market power mitigation mechanisms, this paper proposes a mitigation-aware strategic bidding model and studies the bidding strategies of the market participants under current practice. The proposed bidding model has a bilevel structure with strategic participant’s profit maximization problem in the upper level and the dispatch problem for market operators in the lower level. In particular, the consideration of potential offer mitigation is incorporated as upper-level constraints based on the conduct and impact tests. This bilevel problem is reduced to a single-level mixed-integer linear program using the KKT optimality conditions,  duality theory, and linearization. 
Numerical results illustrate how a strategic player can exercise market power to achieve a higher profit even under the current market power mitigation process and we analyze the social impact that the market power exercise results.
\end{abstract}

\begin{IEEEkeywords}
Bidding strategy, electricity market, market power mitigation.
\end{IEEEkeywords}

\section{Introduction}
The restructuring of the traditional monopoly-based power industry dates back to the early 1980s with the aim of introducing fair competition and improving economic efficiency \cite{David2000a}. 
A typical liberalized market is  hierarchical with a market operator and a group of market participants, \textit{e.g.}, generation companies (GenCos), large consumers, and renewable investors.
In theory, well-defined markets can lead to perfect competition and maximized social welfare. 
However,  existing  electricity markets have shown their vulnerability to price distortion and market manipulation.
The ability of a single market participant (or group of participants) to influence price and distort the market  is referred to as \emph{market power} \cite{Kirschen2004a}.

In the wholesale electricity market, as a profit-seeking entity, a strategic participant may exercise market power via two different approaches: \emph{economic withholding} and \emph{physical withholding}. 
Economic withholding involves submitting strategic bids that deviate from the true marginal cost or utility. It has been observed that  GenCos \cite{Ruiz2009, Pozo2011} and large consumers \cite{Kazempour2015b} frequently employ economic withholding. 
Physical withholding entails a participant \emph{not} offering to sell or \emph{not} scheduling an output according to its actual capacity or load. For example, a GenCo can intentionally create electricity scarcity and strategically curtail renewable output \cite{AzizanRuhi2018}, and energy storage devices can optimize arbitrage opportunities \cite{Wang2017}. 
According to the market reports in \cite{PotomacEconomics2022a, PotomacEconomics2022}, there can be observed market-wide potential economic and physical withholding at \(\sim \)2\% and  \(\sim \)3\% of capacity, respectively, and a noticeable growing trend over the years in the NYISO and MISO.
This work focuses attention on economic withholding.

As an attempt to hedge risks and maximize profit, market participants often embed predicted clearing results into their decision-making process. 
A bilevel optimization problem is commonly used to build such a bidding strategy \cite{Ruiz2009, Pozo2011, Dai2017, Wang2017, AzizanRuhi2018, Wang2021b}.
The upper-level (UL) problem maximizes the profit of the strategic player calculated with the generation dispatch and clearing prices from the market clearing problem in the lower level. 
The lower-level (LL) problem determines the market outcome based on the offer from the strategic player as well as those from its competitors.

Market power mitigation has been a persistent challenge in designing liberalized wholesale electricity markets. Significant effort has been made to investigate the potential for strategic bidding and the design of specialized mitigation mechanisms. 
In \cite{Wang2021b}, a penalty charge is introduced to redistribute the increased payments from the newly occurring congestion and ease its negative impact.
In \cite{Guo2019a}, a novel market power mitigation clearing mechanism based on pre-determined bidding capacity division is proposed to limit potential market power execution.
However, the feasibility of implementation and effectiveness of these methods have not been thoroughly investigated. 
On the other hand, among the existing electricity markets in the US, two fundamental approaches are utilized to mitigate market power: the \emph{structural approach} and the \emph{conduct and impact (C\&I) approach}. The former approach checks for the existence of pivotal players according to their ability to relieve congestion along certain transmission lines. For the latter case, a two-step process is employed first to detect a player's strategic conduct and then to assess the impact of such conduct on the market-clearing price. 
More details can be found in \cite{Graf2021}.
The performance of these two fundamentally different approaches is still a topic of intense debate \cite{Bushnell2019}. This paper  examines the C\&I approach.

Very few studies have paid attention to the effects of market power mitigation from a strategic bidding perspective. Based on this observation, we propose a mitigation-aware bidding model to investigate the feasibility of exercising market power under the existing market power mitigation mechanisms. By studying strategic bidding behavior, we seek to obtain insight into the influence and effectiveness of existing mitigation approaches. GenCos are considered as the strategic entities in this paper. 
The proposed mitigation-aware bidding model is formulated as a bilevel optimization problem. Specifically, in the upper level, the bidding strategy of a strategic participant is developed accounting for the potential of offer mitigation. The bilevel model is later converted into a single-level equivalent and linearized with the KKT optimality conditions, duality theory and well-known linearization schemes. 

The main contributions of this paper are two-fold:
\begin{enumerate}
\item The development of a simple optimization-based model that demonstrates the driving forces of exercising market power and the profitability of a strategic market participant under offer mitigation processes. The mitigation-aware strategic player can gain additional profit by taking advantage of its market share as well as network congestion. Our model also captures the fact that non-strategic GenCos may benefit   in the presence of  mitigation-aware strategic market participants.
\item Illustration of the vulnerability of electricity markets to market power exercise with limited offer mitigation tools. Even non-strategic players have insufficient incentive to resist the exercise of market power. The proposed mitigation-aware bidding framework can serve as an analysis tool for alternative market designs.
\end{enumerate}

\section{Model Formulation}\label{sec:model}
\subsection{Mitigation-Unaware Bidding Strategy}
In the electricity market, the market operator solves a dispatch problem to maximize social welfare (or minimize the total generation cost in case of inelastic demand).
Consider a power network modeled as a graph $\mathcal{G}:=(\mathcal{N},\mathcal{E})$, where each edge $(m,n) \in \mathcal{E}$ represents a branch, and each node $m\in \mathcal{N}$ represents a bus.  For each branch $(m,n)\in \mathcal{E}$, let $p_{mn}$ denotes the power flow from bus $m$ to $n$.
For each bus $m$, let $D_m$ denote the aggregate load and $\theta_m$ the voltage phase angle.
We assume there are $N$ GenCos in the market. However, we take the perspective of considering a single GenCo $G$, and its associated set of generating units $\Omega_{G}$. The remaining $N-1$ GenCos are considered as a single entity with generating units in the set $\Omega'_{G}$. For each unit $i\in \Omega_G\cup \Omega'_G$,  $g_i$  denotes the generation output. The generation dispatch problem of the market operator is formulated using the DC power flow model as follows:
\allowdisplaybreaks
\begin{subequations}
      \label{eq:dcopf}
        \begin{align}
       \min_{\substack{\Xi^{LL}}}
        & \quad   \sum_{i\in \Omega_{G}}  \hat{c}_{i}g_{i}+\sum_{j\in \Omega'_{G}}  \hat{c}_{j}g_{j}       \label{eq:opfobj}  \\
        \mathrm{s.t.}& \quad\sum_{i\in (\Omega_{G} \cap I_m)}g_{i}+\sum_{j\in (\Omega'_{G} \cap I_m)}g_{j} = D_{m} + \sum_{n:m\rightarrow n}p_{mn} \nonumber\\
        &  \qquad \quad\quad\quad\quad\quad\quad\quad-\sum_{l:l\rightarrow m}p_{lm}: \lambda_{m},~~\forall m \in\mathcal{N}  \label{eq:balance}\\
       & \quad p_{mn}=B_{mn}(\theta_{m}-\theta_{n}) :  \nu_{mn}, \quad \forall (m,n) \in\mathcal{E}   \label{eq:dcflow} \\
       &  \quad -\overline{P}_{mn}\le p_{mn} \le \overline{P}_{mn}: \sigma^-_{mn}, \sigma^+_{mn}, ~~ \forall (m,n) \in\mathcal{E}   \label{eq:dcbound}\\
       &  \quad 0\le g_{i} \le \overline{G}_{i}: \mu^-_{i}, \mu^+_{i}, \quad \forall i\in \Omega_{G} \label{eq:gbound1}\\
       &  \quad 0\le g_{j} \le \overline{G}_{j}: \mu'^-_{j}, \mu'^+_{j}, \quad \forall j\in \Omega'_{G} \label{eq:gbound2}\\
       &  \quad -\pi \le \theta _{m} \le \pi : \delta^-_{m}, \delta^+_{m}, \quad \forall m \in\mathcal{N} \label{eq:thetabound}
      \end{align} 
\end{subequations}
where $\Xi^{LL} :=\{g_{i},g_{j},p_{mn},\theta_{m} \}$    is the set of LL decision variables.  Note that generation variables have been partitioned into  sets corresponding to GenCo $G$ and the lumped GenCos (this will facilitate future analysis). Equation~\eqref{eq:opfobj} minimizes the total generation cost, \eqref{eq:balance} represents the nodal supply and demand balance,
\eqref{eq:dcflow} is the linear approximation of the line flow, \eqref{eq:dcbound} enforces the transmission capacity limits of each line, \eqref{eq:gbound1} and \eqref{eq:gbound2} are generation bounds for units, \eqref{eq:thetabound} imposes voltage angle bounds for each bus. 
$I_m$ identifies the generating units connected to bus $m$. The associated dual variables are indicated  following the respective equations. In particular, $\lambda_{m}$, the dual variable associated with \eqref{eq:balance}, is interpreted as the market clearing price, \textit {i.e.}, locational marginal price (LMP).  Generally, a market clearing price cap $\overline{\lambda}$ is necessary to guarantee the availability and accessibility of electricity under certain extreme conditions \cite{FERC2016}.

From the strategic GenCo's perspective, each GenCo develops its bidding strategies taking the market outcome into account using the bilevel bidding model. Thus, the strategic player maximizes its profit in the upper level and incorporates the market dispatch problem \eqref{eq:dcopf} in the lower level. Market participants generally bid with pairs of the offer price and quantity $(\hat{c}_{i},\hat{g}_{i})$ in the electricity market. 
Note that only economic withholding scenarios are discussed in this paper; thus, it is assumed that the quantity offers are submitted as nominal values, \textit{i.e.}, $\hat{g}_{i} = \overline{G}_{i},~\forall i \in \Omega_G$. 
The bilevel problem of the strategic GenCo $G$ is formulated as follows:
\begin{subequations}
      \label{eq:bid}
        \begin{align}
          \max_{\substack{\hat{c}_{i}}, \lambda_m, g_i}    & \quad   \sum_{i\in \Omega_{G}}  \left(\lambda_{m(i)}-c_{i}\right)g_{i}   \label{eq:ulobj}\\
          \mathrm{s.t.} 
           & \quad 0\le \hat c_{i} \le \overline{c}, \quad\forall i\in \Omega_{G}\label{eq:offercap}\\
& \quad 0\le \lambda_{m} \le \overline{\lambda}, \quad\forall m\in \mathcal{N} \label{eq:pricecap}\\           
        & \quad \lambda_{m}, g_{i} \in \textrm{arg}                \min_{\substack{\Xi^{LL}}}
        \quad   \sum_{i\in \Omega_{G}}  \hat{c}_{i}g_{i}+\sum_{j\in \Omega'_{G}}  \hat{c}_{j}g_{j}       \label{eq:llobj}  \\
        &  \quad    \quad          \quad    \quad\quad\quad\mathrm{s.t.}
       \quad   \eqref{eq:balance}-\eqref{eq:thetabound}
      \end{align} 
\end{subequations}
where the objective function \eqref{eq:ulobj} represents the profit of the considered GenCo, $\lambda_{m(i)}$ is the clearing price of unit $i$ under bus $m$, and \eqref{eq:offercap} and \eqref{eq:pricecap} represents the market offer and clearing price cap, respectively.

\subsection{Mitigation-Aware Bidding Strategy}

We  incorporate the consideration of possible mitigation policies such  as C\&I assessment into a strategic player's optimal bidding model.
The strategic bidding model \eqref{eq:bid} presented above is referred to as the mitigation-unaware bidding model, and the one accounting possible offer mitigation is the mitigation-aware bidding model.

\begin{figure}[t]
        \centering
        \includegraphics[width=0.9\columnwidth]{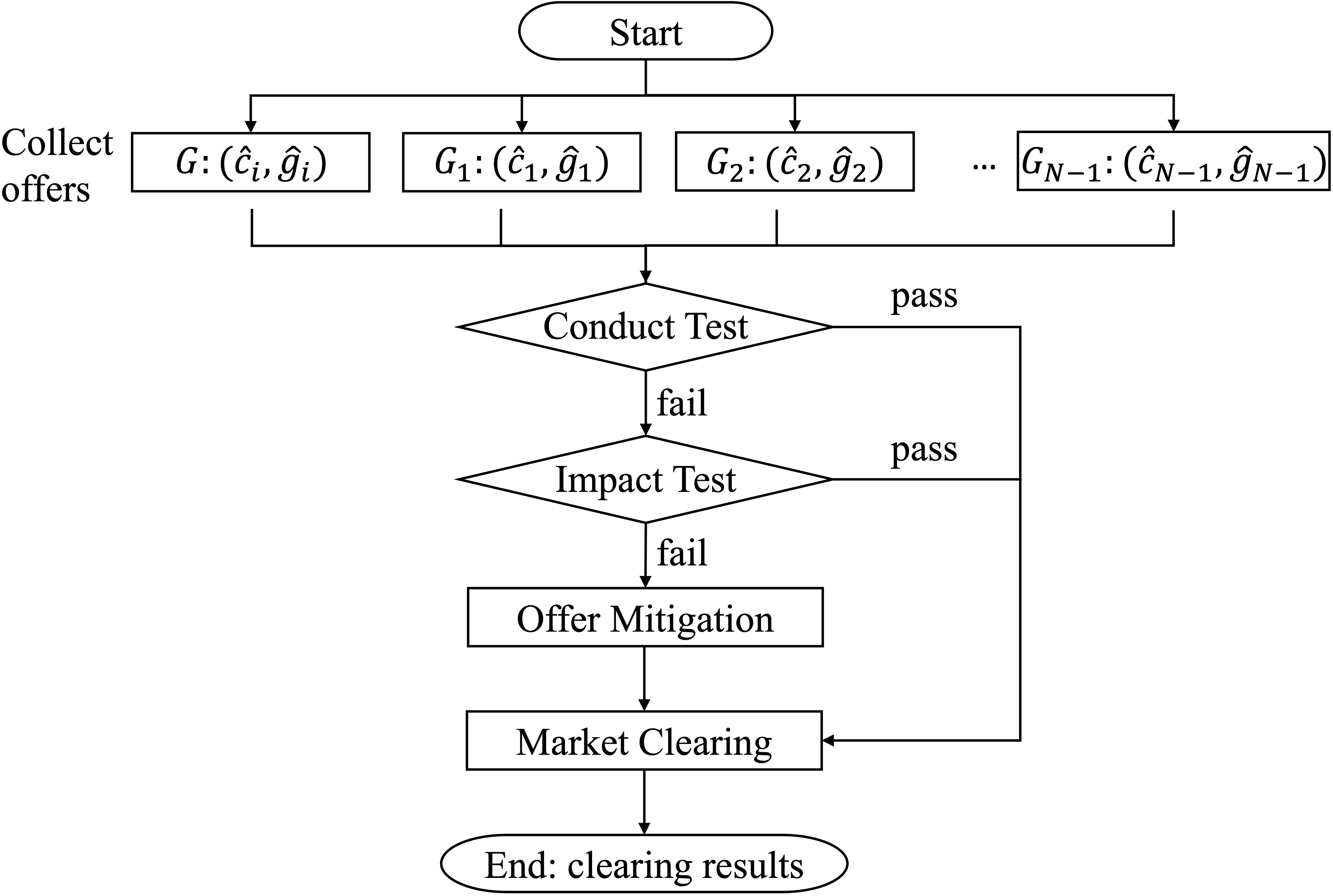}
        \vspace{-2mm}
        \caption{Market clearing workflow with market power mitigation process based on the C\&I tests.  GenCos $G_1- G_{N-1}$ correspond to the lumped GenCos.}
        \vspace{-3mm}
        \label{fig:workflow}
\end{figure}

Fig.~\ref{fig:workflow} gives the workflow for one market clearing period.
At the beginning of the clearing period, the market operator collects offers from the market participants and performs a two-step C\&I assessment. In the \emph{conduct test}, participants' offers are compared to the reference levels maintained by the market operator. If the submitted offer exceeds a specific threshold, it fails the conduct test. In this case, an \emph{impact test} will be executed. This test evaluates the impact of the conduct-test-failed offers by comparing the resulting market prices to the ones obtained from replacing the submitted offers with their reference levels. Note that the impact is determined collectively instead of unit-specifically \cite{NYISO2021}. If an offer fails both tests, it will be mitigated to its reference level before it is transferred to the final market clearing process. The reference levels used throughout the mitigation process are determined based on units' incremental costs, previously accepted offers, previous market prices or negotiated rates \cite{Graf2021}.
Offer mitigation is typically imposed with different thresholds in different areas: constrained and unconstrained areas. Generally, constrained areas correspond to  in-city or frequently congested regions, thus having higher shadow prices. Constrained areas are typically assigned more stringent thresholds, determined by the annual average prices and constrained hours; for unconstrained areas, a widely-adopted mitigation threshold is 300\% higher than the reference level in the conduct test and 200\% in the impact test.

Clearly it is in the interest of a strategic player to implement a bidding strategy that can circumvent the mitigation process. A generated offer can bypass the mitigation process if it is able to satisfy the following constraints:
\begin{subequations}
        \begin{align}
         & \quad |\hat{c} _{i}-c_{i}^0| \le x_{i}, \quad \forall i \in \Omega_{G} \label{eq:c_ref}\\
      & \quad |\lambda_{m}-\lambda_{m}^0| \le y_{m},\quad \forall m \in \mathcal{N}\label{eq:price_ref}
      \end{align} 
    \end{subequations}
where \eqref{eq:c_ref} and \eqref{eq:price_ref} correspond to the conduct test and the impact test, respectively. $c_{i}^0$ represents the estimated reference level of offer price, and $\lambda_{m}^0$ indicates the estimate for a competitive market clearing price. $x_{i}$ and $y_{m}$ are test thresholds. Though the exact reference levels are not provided to the market participants, strategic participants can make estimates using public information. 

Ideally, including one of the constraints \eqref{eq:c_ref} or \eqref{eq:price_ref} in the UL problem will be sufficient to avoid the possible mitigation. It is worth noting that, as reported in \cite{Patton2018}, only approximately less than one percent of the offers failing the conduct test will eventually trigger the impact test. Hence, in terms of bypassing offer mitigation, applying \eqref{eq:c_ref} can be understood as a more conservative movement. In other words, the strategic player may be in a better position of gaining higher profits if they can accurately estimate the impact of their bidding behavior on the market clearing prices.
Therefore, the mitigation-aware bidding problem can be formulated as follows:
\begin{subequations}
\label{eq:bid-aware}
        \begin{align}
          \max_{\substack{\hat{c}_{i}}, \lambda_m, g_i}    & \quad   \sum_{i\in \Omega_{G}}  \left(\lambda_{m(i)}-c_{i}\right)g_{i}  \label{eq:ulobj_a}\\
          \mathrm{s.t.} 
           & \quad \eqref{eq:offercap},~\eqref{eq:pricecap},~\eqref{eq:c_ref},~\eqref{eq:price_ref} \label{eq:ulcon_a} \\
        & \quad \lambda_{m}, g_{i} \in \textrm{arg}               \min_{\substack{\Xi^{LL}}}
        \quad   \sum_{i\in \Omega_{G}}  \hat{c}_{i}g_{i}+\sum_{j\in \Omega'_{G}}  \hat{c}_{j}g_{j}       \label{eq:llobj_a}  \\
         &  \quad    \quad        \quad\quad\quad\quad\mathrm{s.t.}
       \quad   \eqref{eq:balance}-\eqref{eq:thetabound}. \label{eq:llcon_a}
      \end{align} 
      \vspace{-6mm}
\end{subequations}

\subsection{Single-Level Equivalent}

Bilevel problems are strongly NP-hard, and it has been proven an NP-hard task to merely evaluate the optimality of a solution, even for linear bilevel problems as \eqref{eq:bid} or \eqref{eq:bid-aware} \cite{Hansen1992, Vicente1994}. 
A common approach to handle linear bilevel problems such as \eqref{eq:bid} and \eqref{eq:bid-aware} is to derive and solve  single-level equivalents.
As in \cite{Arroyo2010}, the original bilevel problem \eqref{eq:bid-aware} is first converted to an equivalent single-level problem using the KKT optimality conditions, giving:
\begin{subequations}
 \label{eq:mpec}
    \begin{align}
          \max_{\substack{\Xi^{SL}}}    & \quad   \sum_{i\in \Omega_{G}}  \left(\lambda_{m(i)}-c_{i}\right)g_{i}  \label{eq:mpecobj}\\
          \mathrm{s.t.} 
          & \quad \eqref{eq:balance},~\eqref{eq:dcflow},~\eqref{eq:offercap},~\eqref{eq:pricecap},~\eqref{eq:c_ref},~\eqref{eq:price_ref}  \\
          & \quad \hat{c}_{i}-\lambda_{m(i)}-\mu_{i}^-+\mu_{i}^+=0 , \quad \forall i \in \Omega_{G}\label{eq:lgg}\\
          & \quad \hat{c}_{j}-\lambda_{m(j)}-\mu'^-_{j}+\mu'^+_{j}=0 , \quad \forall j \in \Omega'_{G}\label{eq:lggc}\\
          & \quad \lambda_{m}-\lambda_{n}-\nu_{mn}-\sigma_{mn}^-+\sigma_{mn}^+=0, \quad \forall(m,n) \in \mathcal{E}\label{eq:lgpf}\\
          & \quad \!\sum_{(m,n) \in\mathcal{E} }(\nu_{mn}B_{mn}-\nu_{nm}B_{nm})-\delta^-_{m}+\delta^+_{m}=0, ~~ \forall m \in \mathcal{N}\label{eq:lgtheta}\\
          & \quad 0\le \sigma_{mn}^-\perp \overline{P}_{mn}+p_{mn} \ge0, \quad \forall(m,n) \in \mathcal{E} \label{eq:cc1}\\
          & \quad 0\le \sigma_{mn}^+\perp p_{mn}-\overline{P}_{mn}  \ge0, \quad \forall(m,n) \in \mathcal{E}\\
          & \quad 0\le \mu_{i}^-\perp g_{i}\ge0, \quad \forall i \in \Omega_{G}\\
          & \quad 0\le \mu_{i}^+\perp  g_{i}-\overline{G}_{i}\ge0, \quad \forall i \in \Omega_{G}\\
          & \quad 0\le \mu'^-_{j}\perp g_{j}\ge0, \quad \forall j \in \Omega'_{G}\\
          & \quad 0\le \mu'^+_{j}\perp  g_{j}-\overline{G}_{j}\ge0, \quad \forall j \in \Omega'_{G}\\
          & \quad 0\le \delta_{m}^-\perp \pi+\theta_{m} \ge0, \quad  \forall m\in \mathcal{N}\\
          & \quad 0\le \delta_{m}^+\perp \theta_{m}-\pi   \ge0, \quad \forall m\in \mathcal{N}\label{eq:cc2}
    \end{align}
\end{subequations}
where $\Xi^{SL} :=\{\hat{c}_i,g_{i},g_{j},p_{mn},\theta_{m} ,\lambda_{m},\nu_{mn},\sigma^-_{mn}, \sigma^+_{mn},\mu^-_{i}, \mu^+_{i},\\ \mu'^-_{j}, \mu'^+_{j},\delta^-_{m}, \delta^+_{m}\}$    is the set of decision variables for~\eqref{eq:mpec}.  Equation~\eqref{eq:lgg}\textendash\eqref{eq:lgtheta} correspond to the dual constraints of the LL problem \eqref{eq:llobj_a} and \eqref{eq:llcon_a}, and constraints \eqref{eq:cc1}\textendash\eqref{eq:cc2} enforce the complementary conditions.
The notation $\perp$ denotes orthogonality in addition to the stated inequalities.

The  single-level equivalent problem \eqref{eq:mpec} is nonconvex. To make \eqref{eq:mpec} tractable, the bilinear terms are linearized through optimality conditions \cite{Ruiz2009} and the Special Ordered Sets of Type 1 (SOS1) variables \cite{Siddiqui2013}. The linearization techniques are explained in Appendix~\ref{secA:linear}.

\section{Numerical Experiments}\label{sec:sim}

\begin{figure}[t]
    \centering
    \includegraphics[width=0.45\linewidth]{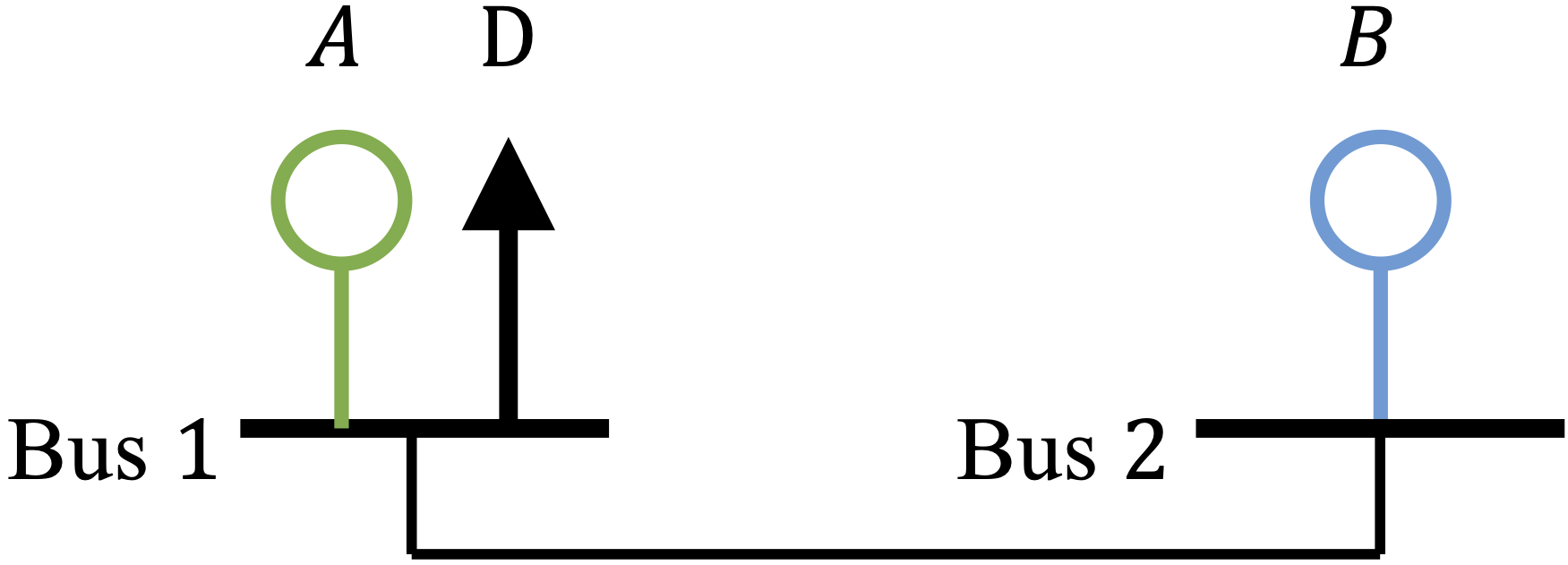}
    \vspace{-2mm}
    \caption{2-bus test system.}
    \vspace{-3mm}
    \label{fig:2bus}
\end{figure}

Our case study is performed on a 2-bus test system, as depicted in Fig.~\ref{fig:2bus}. Due to lack of space, more involved examples are relegated to the Appendix.
Four bidding strategies are tested and compared:
\begin{enumerate}
    \item \textbf{Non-strategic bidding:} Market participants submit offers with true marginal costs.
    \item \textbf{Mitigation-unaware bidding:} Market participants submit offers based on \eqref{eq:bid}.
    \item \textbf{Conduct-aware bidding:} Market participants submit offers based on \eqref{eq:bid-aware} with only \eqref{eq:c_ref}.
    \item \textbf{Impact-aware bidding:} Market participants submit offers based on \eqref{eq:bid-aware} with only \eqref{eq:price_ref}.
\end{enumerate}

The case with two GenCos is considered: Unit  $A$ is the strategic unit, and Unit $B$ is the competitor. There is at most one strategic unit and one load at each bus. 
The marginal cost $c_i$ for each unit is \$20/MWh and generation capacity $\overline{G}_i$ is 30 MW.
Demand is set at $D=50$ MW. 

Perfect prediction is assumed as market participants accurately estimate their competitors' bidding behaviors, and the market operator sets the reference levels for offer mitigation at participants' true operational parameters. 
The mitigation thresholds for C\&I tests are assumed to be 100\% higher than the reference levels \cite{Graf2021}. The market offer cap $\overline{c}$ and clearing price cap $\overline{\lambda}$ are set at \$100/MWh and \$200/MWh, respectively. For strategic players, the estimation for reference price $\lambda_{m}^0$ is carried out by solving the dispatch problem with the reference level offers. In practice,   to guarantee fair dispatch among the price-tied units,
the market operator applies ``tie-breaking'' constraints to price-tied units. We include these in our model; for clarity of exposition, the details are deferred to Appendix~\ref{secA:tie}.

\subsection{Mitigation-unaware bidding \& market power mitigation}

We first consider various bidding strategies for the strategic unit $A$. The competitor $B$ is assumed to be non-strategic.
The bidding, mitigation and clearing results are given in Table~\ref{table:2bus_mit}.  
``Profit$^*_i$'' suggests the expected profits for the strategic units as they make the bidding decisions. Dashes in the table mean that the mitigation is not triggered and the clearing results remain unchanged after mitigation.

 For the case with non-strategic bidding, the clearing prices are the same as the submitted offers and units' marginal costs; hence, there is no explicit profit through market clearing. Given this context, the strategic unit has the incentive either not to bid in the market or raise its offer price to increase the clearing price. 
Assuming that a strategic unit always attempts to expand its profits, it may seek to adopt the mitigation-unaware bidding model and exercise economic withholding. As a result, it can be seen that Unit $A$ offers at the maximum acceptable price $\overline{c}$, raises the clearing price and makes a high profit. Meanwhile, the non-strategic unit also benefits from Unit A's strategic behavior. Indeed, this somewhat counterintuitive outcome is observed in practice \cite{Borenstein2000}.

\begin{table} [t]
        \caption{Clearing Results in the Uncongested Network\rlap{\textsuperscript{1}}} \label{table:2bus_mit} 
        \vspace{-2mm}
        \centering 
        \setlength\tabcolsep{4pt}
        {
        \begin{threeparttable}
                \begin{tabular}{cccccccccc}
                        \hline
\multirow{2}{*}{\begin{tabular}[c]{@{}c@{}}Strategy of \\ Unit $A$\end{tabular}}                                                      & \multirow{2}{*}{Unit} & \multicolumn{4}{c}{Before Mitigation}                                                                                                                                                                                                      & \multicolumn{4}{c}{After Mitigation}                                                                                                                                                      \\ \cmidrule(r){3-6} \cmidrule(r){7-10}
                         & & $\hat{c}_i$& $g_i$& $\lambda_i$ & Profit$^*_i$&$\hat{c}_i$ &$g_i$ &$\lambda_i$ & Profit$_i$ \\ \hline
                        \multirow{2}{*}{Non-Strategic}                                                    & $A$    & 20      & 25     & 20      & 0        & -     & -       & -      & 0       \\
                                                                          & $B$    & 20      & 25     & 20      & 0        & -       & -       & -      & 0       \\
\multirow{2}{*}{\begin{tabular}[c]{@{}c@{}}Mitigation-\\ Unaware\end{tabular}}                                                    & $A$    & 100     & 20     & 100     & 1600     & 20      & 25      & 20     & 0       \\
                                                                          & $B$    & 20      & 30     & 100     & 2400     & -       & 25      & 20     & 0       \\
\multirow{2}{*}{Conduct-Aware} & $A$    & 40      & 20     & 40      & 400      & -       & -       & -      & 400     \\
                                                                          & $B$    & 20      & 30     & 40      & 600      & -       & -       & -      & 600     \\
\multirow{2}{*}{Impact-Aware}  & $A$    & 40      & 20     & 40      & 400      & -       & -       & -      & 400     \\
                                                                          & $B$    & 20      & 30     & 40      & 600      & -       & -       & -      & 600\\ \hline
                        \end{tabular}                        
                        \begin{tablenotes}

\item[1]Units: $\hat{c}_i$ {[}\$/MWh{]}, $g_i$ {[}MW{]}, $\lambda_i$ {[}\$/MWh{]}, Profit$_i$ {[}\${]}.

\end{tablenotes}
                        \end{threeparttable}
        }
        \vspace{-3mm}
\end{table}

Strategic behaviors are closely monitored by the market operator. We apply the mitigation process to the submitted offers, and the results are summarized in the right half of Table~\ref{table:2bus_mit}. 
It can be seen that the strategic offers generated from the mitigation-unaware bidding model are vulnerable to the offer mitigation process. 
Once detected in the mitigation process, bidding offers are reset to the reference level, which leaves the strategic unit with a relatively narrow profit margin.
This motivates market participants to adopt a smarter bidding model.
\subsection{Effects of mitigation-aware bidding}

The clearing results adopting mitigation-aware bidding strategies are shown in Table~\ref{table:2bus_mit}. It can be seen that mitigation-aware bidding (either conduct-aware or impact-aware strategies) can successfully bypass the offer mitigation and achieve a higher profit. 
Compared to the mitigation-unaware bidding, the conduct-aware bidding offers at \$40/MWh instead of \$100/MWh and yields \$400 higher profit for Unit $A$. The clearing results from the conduct- and impact-aware bidding are identical due to the same operational parameters for units.

\begin{table}[t]
        \caption{Post-Mitigation Clearing Results in the Uncongested Network} \label{table:2bus_ci} 
        \vspace{-2mm}
        \centering 
        \setlength\tabcolsep{11pt}
        {
        \begin{threeparttable}
                \begin{tabular}{cccccc}
                        \hline
                        \begin{tabular}[c]{@{}c@{}}Strategy of \\ Unit $A$\end{tabular}& Unit& $\hat{c}_i$ & $g_i$ & $\lambda_i$ & Profit$_i$  \\ \hline
\multirow{2}{*}{Conduct-Aware} & $A$ & 20-$\epsilon$          & 30    & 20          & 300     \\
 & $B$ & 20          & 20    & 20          & 0     \\
\multirow{2}{*}{Impact-Aware} & $A$ & 40          & 20    & 40          & 600     \\
& $B$& 20          & 30    & 40          & 600 \\ \hline
                        \end{tabular}
                        \end{threeparttable}
                        }
                        \vspace{-3mm}
\end{table}

What is worth noting is the lower profits of Unit $A$ compared to that of Unit $B$ at the end of the clearing period. The reason for this  is that when Unit $A$ submits an offer higher than its true marginal cost, Unit $B$ becomes the first-to-clear unit; then, Unit $A$ is dispatched with 20 MW and Unit $B$ gets 30 MW. As a consequence, with the same clearing prices, Unit $B$ makes higher profits even after getting mitigated and cleared with a reference-level offer, \textit{i.e.}, true marginal cost level. This result is also somewhat counterintuitive. However, the relatively lower profits still promise a better return compared to the outcome from non-strategic or mitigation-unaware bidding strategies, yielding the highest profit among all bidding strategies. What's more, if the competitor also adopts a mitigation-aware bidding strategy, two units will evenly supply  the total demand and make a profit of \$500, respectively. In other words, mitigation-aware bidding strategies lead to a suboptimal outcome. It also shows that the capacity limit is one major source of exercising market power.

\subsection{Comparison between conduct- \& impact-aware bidding strategies}

The difference between the conduct-aware and impact-aware bidding strategies becomes clearer when the heterogeneous generating units are considered, \textit{i.e.}, units have different marginal costs. 
We now set the marginal cost of Unit $A$ at \$10/MWh and the competitor unit $B$ at \$20/MWh. The results are shown in Table~\ref{table:2bus_ci}.
Under the conduct-aware bidding strategy, Unit $A$ obtains its offer at \$(20-$\epsilon$)MWh (Note, $\epsilon>0$ is a small value that causes the offer price to be slightly lower than the nominal value. Although we don't specify this everywhere, it should be noted that  this holds for all offers with mitigation-aware bidding strategies.) and gets cleared with 30 MW. Accounting for the impact test, Unit $A$ submits its offer at \$40/MWh and closes the deal with 20 MW. Though the final generation output is lower, impact-aware bidding ends up with higher profits. The conclusion can be drawn that when the strategic unit is not at the marginal position according to the true marginal cost, it's rational to bid under the impact-aware strategy and pursue a higher profit.
As such, the conduct-aware bidding strategy is shown to be a more conservative strategy than the impact-aware strategy. This result is aligned with the fact that only a small portion of the conduct-test-failed offers are eventually mitigated after the impact test in reality. 

\begin{figure}[t]
\centerline{\subfigure[Uncongested Network]{\includegraphics[width=1.7in]{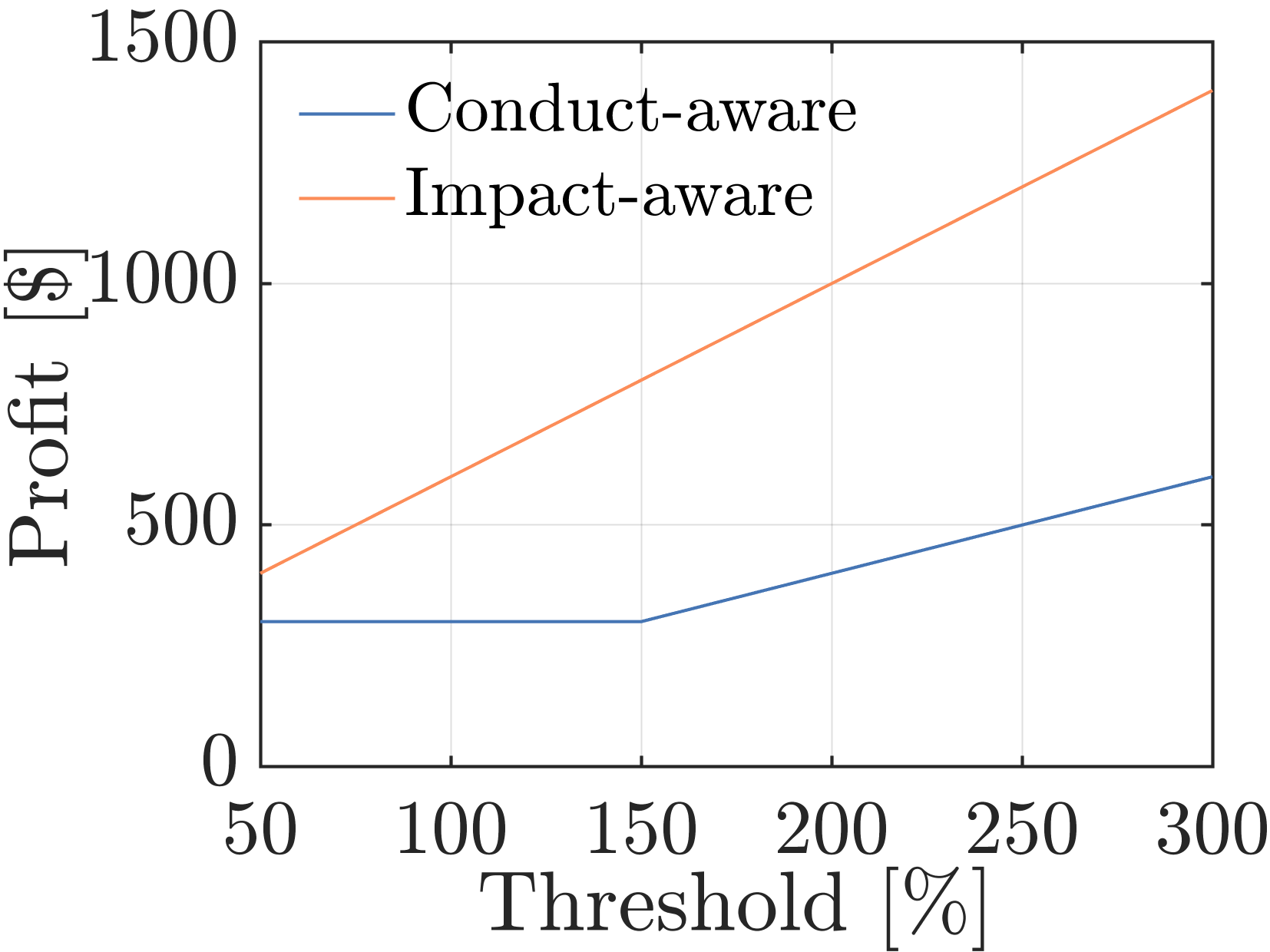}
\label{fig:uncong}
} 
\hfil 
\subfigure[Congested Network]{\includegraphics[width=1.7in]{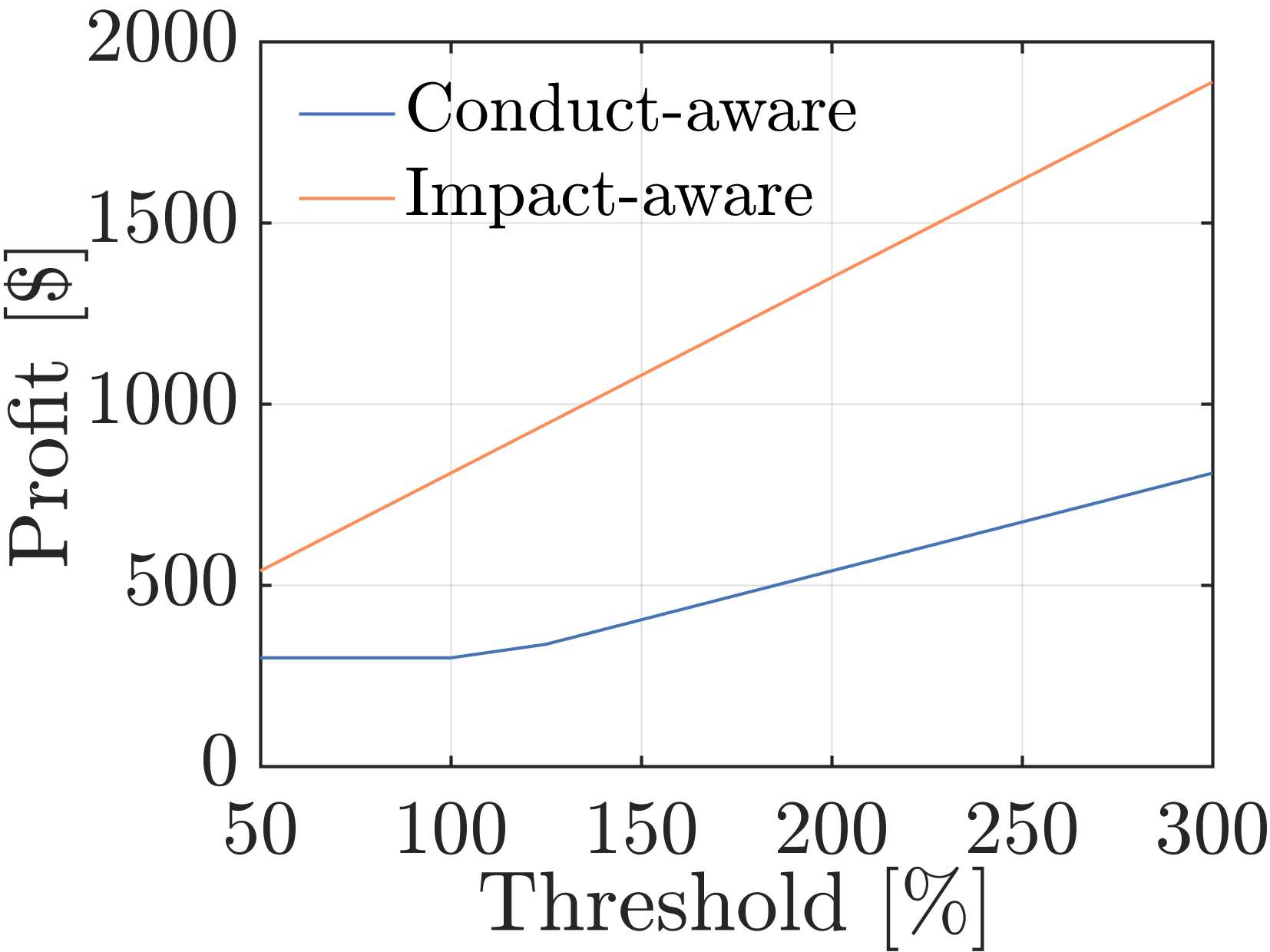}
\label{fig:cong}
}
} 
\vspace{-2mm}
\caption{Post-mitigation clearing results with different mitigation thresholds.}
\label{fig:mit_th} 
\vspace{-3mm}
\end{figure}

We further examine the effects of mitigation thresholds on mitigation-aware bidding. It is easy to imagine that the profits will linearly increase with the growth of thresholds when $c_A=c_B$ or $c_A>c_B$ and be identical in adopting conduct- or impact-aware bidding strategies. The results when $c_A=\$$10/MWh and $c_B=\$$20/MWh are shown in Fig.~\ref{fig:uncong}. The thresholds are selected as 50\%\textendash300\% higher than the reference levels. For the case of conduct-aware bidding, the profit of Unit $A$ remains steady at \$300 when the threshold is relatively lower because its offer price is lower than that of its competitor and gets cleared with higher output at the competitor's offer price. It can be seen that the profit difference between conduct-aware and impact-aware bidding becomes larger as the thresholds become higher. 

On the other hand, market power exercising will eventually jeopardize the social welfare due to the increase in the total generation cost. Consider the case $c_A<c_B$. Once the generating unit $A$ falsely claims a marginal cost $\hat{c}_A>c_B$ when submitting the offer, the market operator will prioritize clearing Unit $B$. Then, the total generation cost is increased compared to the situation of perfect competition. The impact-aware bidding case in Table~\ref{table:2bus_ci} is one typical example of this. Assuming that the marginal utility of demand is \$25/MWh \cite{Kazempour2015b}, then the original social welfare is supposed to be \$550 while it drops to \$450 under mitigation-aware bidding.
Therefore, from the market operator's point of view, it is important to consider the limitation of their mitigation mechanism.

\subsection{Effects of line congestion}
Scenarios without line congestion can be understood as the local competition isolated by exterior congestion. 
We now examine the effects of line congestion. The transmission line limit is set to 23 MW, a value slightly lower than the evenly dispatching decisions. 
For this new scenario, the clearing results are presented in Table~\ref{table:2bus_cong}. First, the difference from the uncongested scenario is reflected in the generation output. Here, the market share remains at 27 MW due to the line congestion. That is, there are chances of high supply dependence in a congestion-rendered isolated local market. Secondly, line congestion causes different clearing prices at different nodes of the congested line. A higher price under this condition can be interpreted as a reward for relieving the congestion. As a result, the suboptimal situation for the strategic unit is resolved. Altogether, the profits for the strategic unit in the congested scenario are higher than those in the uncongested area. 
Hence, congestion serves as a persistent second major source of market manipulation and puts the strategic unit in a better position for profit-seeking.

\begin{table}[t]
        \caption{Clearing Results in the Congested Network} \label{table:2bus_cong} 
        \vspace{-2mm}
        \centering 
        \setlength\tabcolsep{4pt}{
                \begin{tabular}{cccccccccc}
                                                \hline
\multirow{2}{*}{\begin{tabular}[c]{@{}c@{}}Strategy of \\ Unit $A$\end{tabular}}                                                      & \multirow{2}{*}{Unit} & \multicolumn{4}{c}{Before Mitigation}                                                                                                                                                                                                      & \multicolumn{4}{c}{After Mitigation}                                                                                                                                                      \\ \cmidrule(r){3-6} \cmidrule(r){7-10}
                         & & $\hat{c}_i$& $g_i$& $\lambda_i$ & Profit$^*_i$&$\hat{c}_i$ &$g_i$ &$\lambda_i$ & Profit$_i$ \\ \hline
                        \multirow{2}{*}{Non-Strategic}                                                   & $A$    & 20      & 27     & 20      & 0        & -       & -       & -      & 0       \\
                                                                          & $B$    & 20      & 23     & 20      & 0        & -       & -       & -      & 0       \\
\multirow{2}{*}{\begin{tabular}[c]{@{}c@{}}Mitigation-\\ Unaware\end{tabular}}                                                    & $A$    & 100     & 27     & 100     & 2160     & 20      & 27      & 20     & 0       \\
                                                                          & $B$    & 20      & 23     & 20     & 0     & -       & -      & -     & 0       \\
\multirow{2}{*}{Conduct-Aware} & $A$    & 40      & 27     & 40      & 540      & -       & -       & -      & 540     \\
                                                                          & $B$    & 20      & 23     & 20      & 0      & -       & -       & -      & 0     \\
\multirow{2}{*}{Impact-Aware}  & $A$    & 40      & 27     & 40      & 540      & -       & -       & -      & 540     \\
                                                                          & $B$    & 20      & 23     & 20      & 0      & -       & -       & -      & 0\\ \hline
                        \end{tabular}
        }
        \vspace{-3mm}
\end{table}

Fig.~\ref{fig:cong} gives the clearing results considering different mitigation thresholds when the network is congested. The trend for profit increase is similar to that in the uncongested scenario. The main difference is the higher overall profits discussed above and a better chance to gain profits using conduct-aware bidding with lower thresholds attributed to higher output. It is reasonable to conclude that the strategic unit is in a preferable position to exercise market power in a high-demand and frequently congested area. 
Consequently, as in the case of the existing markets, constrained areas are often assigned with a more restrictive mitigation policy to prevent extreme exploitation of market power.

The effects of line congestion are more subtle when the graph representing the physical structure of the network contains loops. A simplified single loop example is described in Appendix~\ref{secA:3bus}. 
A more realistic bidding scenario accounting for multiple agents and multi-block offer curves is presented in Appendix~\ref{secA:6bus}.

\section{Conclusions}\label{sec:conc}
A mitigation-aware strategic bidding model is proposed to investigate the bidding behavior of the strategic market participants under  existing market power mitigation mechanisms and the effectiveness of these practices investigated.
The strategic bidding model is constructed based on a bilevel optimization framework. The UL problem considers profit maximization of the strategic participant, and the LL problem performs market clearing. In particular, the consideration of potential offer mitigation is incorporated into the UL problem.  Numerical results reveal the vulnerability of the electricity market to market power exercise and consequent loss of social welfare with limited mitigation tools. Strategic participant are shown to achieve a higher profit by taking advantage of market scarcity and network congestion and circumventing offer mitigation. Meanwhile, even the non-strategic participants may  benefit from the exercise of market power within the market. 
Future work will  examine strategic behavior in terms of physical withholding and remove the perfect prediction assumption.

\bibliographystyle{IEEEtran}
\bibliography{library.bib}

\begin{thebibliography}{10}
\providecommand{\url}[1]{#1}
\csname url@samestyle\endcsname
\providecommand{\newblock}{\relax}
\providecommand{\bibinfo}[2]{#2}
\providecommand{\BIBentrySTDinterwordspacing}{\spaceskip=0pt\relax}
\providecommand{\BIBentryALTinterwordstretchfactor}{4}
\providecommand{\BIBentryALTinterwordspacing}{\spaceskip=\fontdimen2\font plus
\BIBentryALTinterwordstretchfactor\fontdimen3\font minus
  \fontdimen4\font\relax}
\providecommand{\BIBforeignlanguage}[2]{{%
\expandafter\ifx\csname l@#1\endcsname\relax
\typeout{** WARNING: IEEEtran.bst: No hyphenation pattern has been}%
\typeout{** loaded for the language `#1'. Using the pattern for}%
\typeout{** the default language instead.}%
\else
\language=\csname l@#1\endcsname
\fi
#2}}
\providecommand{\BIBdecl}{\relax}
\BIBdecl

\bibitem{David2000a}
A.~David and {Fushuan Wen}, ``{Strategic bidding in competitive electricity
  markets: a literature survey},'' in \emph{2000 Power Engineering Society
  Summer Meeting (Cat. No.00CH37134)}, vol.~4.\hskip 1em plus 0.5em minus
  0.4em\relax IEEE, 2000, pp. 2168--2173.

\bibitem{Kirschen2004a}
D.~Kirschen and G.~Strbac, \emph{{Fundamentals of Power System
  Economics}}.\hskip 1em plus 0.5em minus 0.4em\relax Chichester, UK: John
  Wiley {\&} Sons, Ltd, 3 2004.

\bibitem{Ruiz2009}
C.~Ruiz and A.~Conejo, ``{Pool Strategy of a Producer With Endogenous Formation
  of Locational Marginal Prices},'' \emph{IEEE Transactions on Power Systems},
  vol.~24, no.~4, pp. 1855--1866, 11 2009.

\bibitem{Pozo2011}
D.~Pozo and J.~Contreras, ``{Finding Multiple Nash Equilibria in Pool-Based
  Markets: A Stochastic EPEC Approach},'' \emph{IEEE Transactions on Power
  Systems}, vol.~26, no.~3, pp. 1744--1752, 8 2011.

\bibitem{Kazempour2015b}
S.~J. Kazempour, A.~J. Conejo, and C.~Ruiz, ``{Strategic Bidding for a Large
  Consumer},'' \emph{IEEE Transactions on Power Systems}, vol.~30, no.~2, pp.
  848--856, 3 2015.

\bibitem{AzizanRuhi2018}
N.~Azizan~Ruhi, K.~Dvijotham, N.~Chen, and A.~Wierman, ``{Opportunities for
  Price Manipulation by Aggregators in Electricity Markets},'' \emph{IEEE
  Transactions on Smart Grid}, vol.~9, no.~6, pp. 5687--5698, 11 2018.

\bibitem{Wang2017}
Y.~Wang, Y.~Dvorkin, R.~Fern{\'{a}}ndez-Blanco, B.~Xu, T.~Qiu, and D.~S.
  Kirschen, ``{Look-Ahead Bidding Strategy for Energy Storage},'' \emph{IEEE
  Transactions on Sustainable Energy}, vol.~8, no.~3, pp. 1106--1127, 2017.

\bibitem{PotomacEconomics2022a}
\BIBentryALTinterwordspacing
{Potomac Economics}, ``{2021 State of the Market Report for the New York ISO
  Markets},'' 2022. [Online]. Available:
  \url{https://www.nyiso.com/documents/20142/2223763/NYISO-2021-SOM-Full-Report-5-11-2022-final.pdf/5307870c-9b62-1720-1708-6b9c157211bb}
\BIBentrySTDinterwordspacing

\bibitem{PotomacEconomics2022}
\BIBentryALTinterwordspacing
------, ``{2021 State of the Market Report for the MISO ELectricity Markets},''
  2022. [Online]. Available: \url{https://cdn.misoenergy.org/20220622 Markets
  Committee of the BOD Item 04 IMM State of the Market Report625261.pdf}
\BIBentrySTDinterwordspacing

\bibitem{Dai2017}
T.~Dai and W.~Qiao, ``{Finding Equilibria in the Pool-Based Electricity Market
  With Strategic Wind Power Producers and Network Constraints},'' \emph{IEEE
  Transactions on Power Systems}, vol.~32, no.~1, pp. 389--399, 1 2017.

\bibitem{Wang2021b}
S.~Wang, X.~Tan, T.~Liu, and D.~H. Tsang, ``{Aggregation of Demand-Side
  Flexibility in Electricity Markets: Negative Impact Analysis and Mitigation
  Method},'' \emph{IEEE Transactions on Smart Grid}, vol.~12, no.~1, pp.
  774--786, 2021.

\bibitem{Guo2019a}
H.~Guo, Q.~Chen, Q.~Xia, and C.~Kang, ``{Market Power Mitigation Clearing
  Mechanism Based on Constrained Bidding Capacities},'' \emph{IEEE Transactions
  on Power Systems}, vol.~34, no.~6, pp. 4817--4827, 11 2019.

\bibitem{Graf2021}
C.~Graf, E.~La~Pera, F.~Quaglia, and F.~A. Wolak, ``{Market Power Mitigation
  Mechanisms for Wholesale Electricity Markets: Status Quo and Challenges},''
  2021.

\bibitem{Bushnell2019}
\BIBentryALTinterwordspacing
J.~Bushnell, S.~M. Harvey, and B.~F. Hobbs, ``{Opinion on System Market Power
  Mitigation},'' Market Surveillance Committee of the California ISO, Tech.
  Rep., 2019. [Online]. Available:
  \url{http://www.caiso.com/documents/msc-draftopiniononsystemmarketpowermitigation-nov5-2019.pdf}
\BIBentrySTDinterwordspacing

\bibitem{FERC2016}
\BIBentryALTinterwordspacing
{FERC}, ``{Offer Caps in Markets Operated by Regional Transmission
  Organizations and Independent System Operators},'' 2016. [Online]. Available:
  \url{https://www.federalregister.gov/documents/2016/12/05/2016-28320/offer-caps-in-markets-operated-by-regional-transmission-organizations-and-independent-system}
\BIBentrySTDinterwordspacing

\bibitem{NYISO2021}
\BIBentryALTinterwordspacing
{NYISO}, ``{Market Participants User ’ s Guide},'' 2021. [Online]. Available:
  \url{https://www.nyiso.com/manuals-tech-bulletins-user-guides}
\BIBentrySTDinterwordspacing

\bibitem{Patton2018}
\BIBentryALTinterwordspacing
D.~Patton, ``{Market Power Mitigation in Wholesale Electricity Markets},''
  2018. [Online]. Available:
  \url{https://www.westerneim.com/Documents/Presentation-MarketPowerMitigationPatton.pdf}
\BIBentrySTDinterwordspacing

\bibitem{Hansen1992}
P.~Hansen, B.~Jaumard, and G.~Savard, ``{New Branch-and-Bound Rules for Linear
  Bilevel Programming},'' \emph{SIAM Journal on Scientific and Statistical
  Computing}, vol.~13, no.~5, pp. 1194--1217, 9 1992.

\bibitem{Vicente1994}
L.~Vicente, G.~Savard, and J.~J{\'{u}}dice, ``{Descent approaches for quadratic
  bilevel programming},'' \emph{Journal of Optimization Theory and
  Applications}, vol.~81, no.~2, pp. 379--399, 5 1994.

\bibitem{Arroyo2010}
J.~Arroyo, ``{Bilevel programming applied to power system vulnerability
  analysis under multiple contingencies},'' \emph{IET Generation, Transmission
  {\&} Distribution}, vol.~4, no.~2, p. 178, 2010.

\bibitem{Siddiqui2013}
S.~Siddiqui and S.~A. Gabriel, ``{An SOS1-Based Approach for Solving MPECs with
  a Natural Gas Market Application},'' \emph{Networks and Spatial Economics},
  vol.~13, no.~2, pp. 205--227, 6 2013.

\bibitem{Borenstein2000}
S.~Borenstein, ``{Understanding Competitive Pricing and Market Power in
  Wholesale Electricity Markets},'' \emph{The Electricity Journal}, vol.~13,
  no.~6, pp. 49--57, 7 2000.

\bibitem{Fortuny-amat1981}
J.~Fortuny-Amat and B.~McCarl, ``{A Representation and Economic Interpretation
  of a Two-Level Programming Problem},'' \emph{The Journal of the Operational
  Research Society}, vol.~32, no.~9, p. 783, 9 1981.

\bibitem{Hu2008}
J.~Hu, J.~E. Mitchell, J.-S. Pang, K.~P. Bennett, and G.~Kunapuli, ``{On the
  Global Solution of Linear Programs with Linear Complementarity
  Constraints},'' \emph{SIAM Journal on Optimization}, vol.~19, no.~1, pp.
  445--471, 1 2008.

\bibitem{Gabriel2010}
S.~A. Gabriel and F.~U. Leuthold, ``{Solving discretely-constrained MPEC
  problems with applications in electric power markets},'' \emph{Energy
  Economics}, vol.~32, no.~1, pp. 3--14, 1 2010.

\bibitem{MISO2021}
\BIBentryALTinterwordspacing
{MISO}, ``{Business practices manual - energy and operating reserve markets},''
  2021. [Online]. Available:
  \url{https://www.misoenergy.org/legal/business-practice-manuals/}
\BIBentrySTDinterwordspacing

\bibitem{Kwon2020}
J.~Kwon, Z.~Zhou, T.~Levin, and A.~Botterud, ``{Resource Adequacy in
  Electricity Markets With Renewable Energy},'' \emph{IEEE Transactions on
  Power Systems}, vol.~35, no.~1, pp. 773--781, 1 2020.

\end{thebibliography}

\newpage 

\appendices

\section{Model Linearization} \label{secA:linear}
\renewcommand{\theequation}{\thesection.\arabic{equation}}
\renewcommand\thefigure{\thesection.\arabic{figure}} 
\renewcommand\thetable{\thesection.\arabic{table}} 
\setcounter{equation}{0}

The  single-level equivalent problem \eqref{eq:mpec} contains nonlinearities as products of continuous decision variables. The corresponding linearization techniques used in this paper are summarized as follows:
\begin{enumerate}
        \item Bilinear products of market-clearing prices and power outputs, \textit{i.e.}, $\lambda_{m}g_{i}$, in the objective function \eqref{eq:mpecobj}: These sets of nonlinear products can be linearized without approximations using the KKT conditions and the strong duality theorem \cite{Ruiz2009}.
        
        The strong duality theorem states that the primal and dual objective function values are equal at the optimum if a problem is convex. Since the LL problem satisfies the prerequisite, the strong duality theorem holds as follows:
        \begin{align}
                \sum_{i\in \Omega_{G}}  \hat{c}_{i}g_{i}=&-\sum_{j\in \Omega'_{G}}  \hat{c}_{j}g_{j} + \sum_{m \in\mathcal{N} }  \lambda_{m} D_{m} \nonumber\\
                &- \sum_{i\in \Omega_{G}} \mu_{i}^+\overline{G}_{i}- \sum_{j\in \Omega'_{G}} \mu'^+_{j}\overline{G}_{j} \nonumber\\
                &-\sum_{(m,n) \in\mathcal{E} }\left(\sigma_{mn}^-\overline{P}_{mn}+\sigma_{mn}^+\overline{P}_{mn}\right) \nonumber\\
                &-\sum_{m \in\mathcal{N} }\left(\delta_{m}^-\pi +\delta_{m}^+\pi\right). \label{eq:sd}
        \end{align} 
  
        Using KKT optimality equalities \eqref{eq:lgg} and \eqref{eq:cc2}, the left hand side of \eqref{eq:sd} can be written  as 
        \begin{align}
                \sum_{i\in \Omega_{G}}  \hat{c}_{i}g_{i} =\sum_{i\in \Omega_{G} } \lambda_{m(i)} g_{i} - \sum_{i\in \Omega_{G}} \mu_{i}^+ \overline{G}_{i}. \label{eq:sub4}
        \end{align}

        Substituting \eqref{eq:sub4} in \eqref{eq:sd} renders
        \begin{align}
                \sum_{i\in \Omega_{G} } \lambda_{m(i)} g_{i} =&-\sum_{j\in \Omega'_{G}}  \hat{c}_{j}g_{j} + \sum_{m \in\mathcal{N} }  \lambda_{m} D_{m}- \sum_{j\in \Omega'_{G}} \mu'^+_{j}\overline{G}_{j} \nonumber\\
                &-\sum_{(m,n) \in\mathcal{E} }\left(\sigma_{mn}^-\overline{P}_{mn}+\sigma_{mn}^+\overline{P}_{mn}\right) \nonumber\\
                &-\sum_{m \in\mathcal{N} }\left(\delta_{m}^-\pi +\delta_{m}^+\pi\right) \label{eq:sd_l}
        \end{align} 
      which is an equivalent linear expression of the bilinear term in the objective function \eqref{eq:mpecobj}.
      
        \item Products of Lagrange multipliers and LL decision variables in the complementarity constraints \eqref{eq:cc1}\textendash\eqref{eq:cc2}: One of the well-known solutions is a Big-M-based reformulation as proposed in \cite{Fortuny-amat1981}. However, difficulties arise when an attempt is made to determine an appropriate $M$ value \cite{Hu2008,Gabriel2010}. An alternative solution is to replace the complementarity constraints using SOS1 variables \cite{Siddiqui2013}. SOS1 variables are defined as a set of variables where at most, one of the variables in the set  can take a non-zero value. Considering the drawback of the Big-M-based method upon implementation, SOS1-based linearization is employed in this paper. 
        Due to lack of space, the complete formulation of the linearized model is neglected here. A simple example of SOS1-based linearization is given below.

Consider a set of complementarity constraints $y^Tg=0$, where $y,g \ge 0$.
According to \cite{Siddiqui2013}, the products $y^Tg=0$ can be recast using a pair of SOS1 variables $v^+, v^-$ constrained as follows: 
        \begin{subequations}
                  \label{eq:sos1}
          \begin{align}
      & u-(v^++v^-)=0\\ 
      & u = \frac {a +b} 2 \\ 
      & v^+-v^- = \frac {a -b} 2 \\ 
      & v^+, v^- \in \textrm{SOS1}.
          \end{align} 
        \end{subequations}
Note that the solution to \eqref{eq:sos1} can lead to a global optimum to the original mathematical programs with complementarity constraints.

\end{enumerate}

After the linearization process described above, a mixed-integer linear program (MILP) is derived from the initial bilevel program and reformulated single-level equivalent. The main advantages include \cite{Arroyo2010}: i) the guaranteed convergence to the optimal solution within a finite number of steps, and ii) the ready availability of well-established branch-and-cut solver.

\section{Tie-Breaking Constraints}\label{secA:tie}
\setcounter{equation}{0}

Tie-breaking constraints ensure that dispatch results are proportional to the submitted capacity. An additional penalty term $\sum_{(i,j)\in \Phi^\mathrm{TB}}\rho \left(\mathrm{tb1} _{ij}+\mathrm{tb2} _{ij}\right)$ is appended to~\eqref{eq:ulobj_a} (where $\Phi^\mathrm{TB}$ represents the set of all price-tied unit candidates, and $\rho$ is a small penalty term). The tie-breaking constraints take the form
 \cite{MISO2021}:
\begin{subequations}
\begin{align}
            &\! \mathrm{tb1}_{ij} \!\ge\! \overline{G}_{i } g_{j}\!-\! \overline{G}_{j} g_{i} , ~~\forall i\!\in\! (\Phi^\mathrm{TB}\!\cap\! \Omega_{G} ), \forall j\! \in\! (\Phi^\mathrm{TB}\!\cap \!\Omega'_{G}) \!\label{eq:tb1}\\
        &\! \mathrm{tb2} _{ij}\!\ge\!  \overline{G}_{j}g_{i} \!-\!\overline{G}_{i}g_{j}, ~~ \forall i\!\in\! (\Phi^\mathrm{TB}\!\cap\!\Omega_{G} ),\forall j \!\in\! (\Phi^\mathrm{TB}\!\cap\! \Omega'_{G})\!\label{eq:tb2}
        \end{align}
\end{subequations}
where $\mathrm{tb1} _{ij}$ and $ \mathrm{tb2} _{ij}$ are non-negative tie-breaking variables. 
Following~\cite{Kwon2020}, we place these constraints in the  UL problem \eqref{eq:ulobj_a} and \eqref{eq:ulcon_a}.

\section{3-Bus Case Study}\label{secA:3bus}
\setcounter{equation}{0}
\setcounter{figure}{0}   
\setcounter{table}{0}

\begin{figure}[t]
    \centering
    \includegraphics[width=0.5\linewidth]{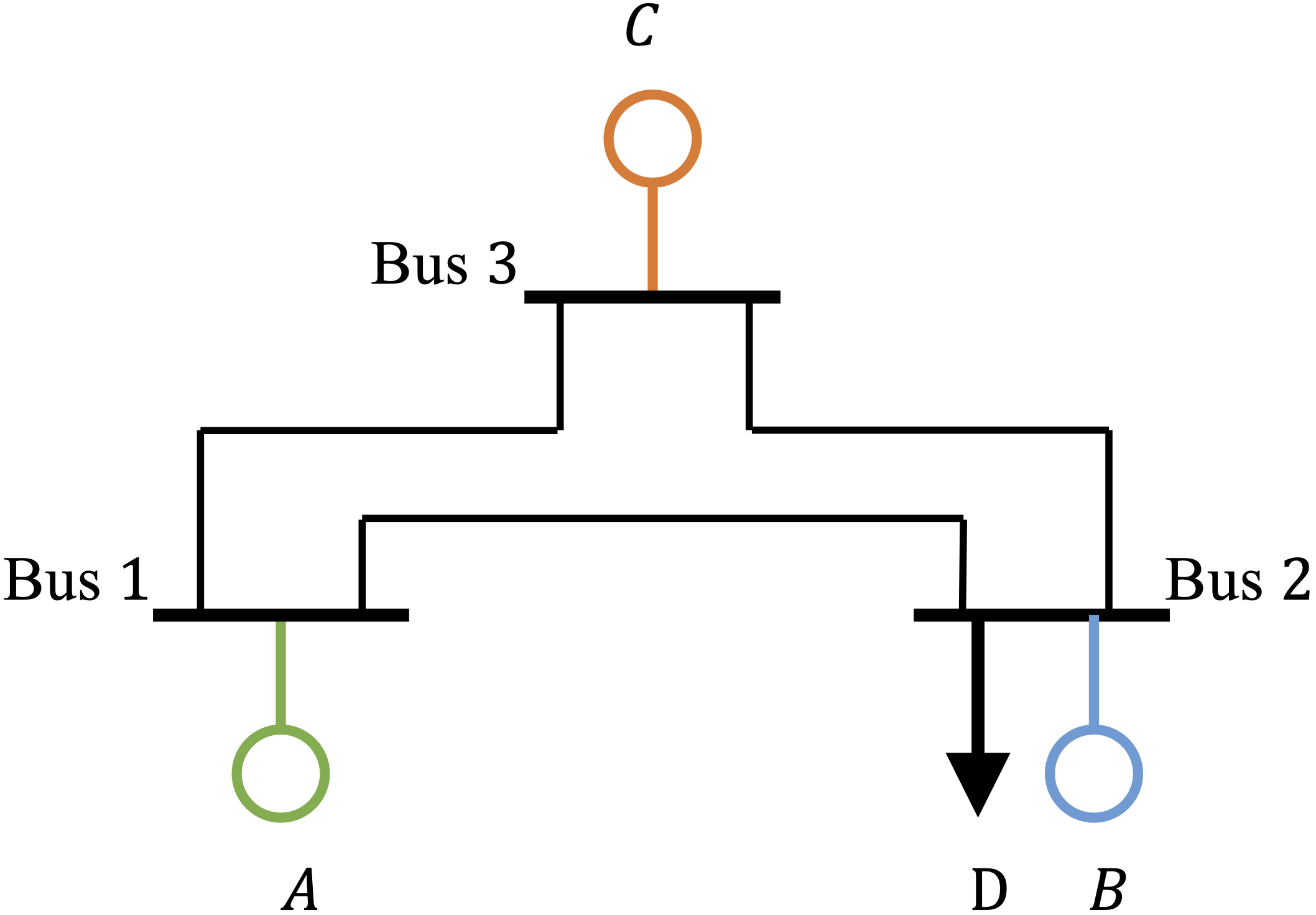}
    \vspace{-2mm}
    \caption{3-bus test system.}
    \label{fig:3bus}
    \vspace{-3mm}
\end{figure}

The effects of congestion on market power exercise are more complicated when there are loops in the network. Take the 3-bus test system as an example, with topology given in Fig.~\ref{fig:3bus}; the other configuration is the same as the 2-bus case study. The clearing results are summarized in Table~\ref{table:3bus_cong}. As Table~\ref{table:3bus_cong} shows, the strategic behavior successfully gets Unit $A$ a profit markup, and, meanwhile, raises the clearing price at $N_2$ and consequently the profit of Unit $B$ in attributions to line congestion. This result exhibits another instance of a suboptimal outcome for the strategic unit.

\begin{table}[t]
        \caption{Clearing Results in the Congested Network with the 3-Bus Test System} \label{table:3bus_cong}
        \vspace{-2mm}
        \centering 
        \setlength\tabcolsep{4pt}{
                \begin{tabular}{cccccccccc}
                        \hline
\multirow{2}{*}{\begin{tabular}[c]{@{}c@{}}Strategy of \\ Unit $A$\end{tabular}}                                                      & \multirow{2}{*}{Unit} & \multicolumn{4}{c}{Before Mitigation}                                                                                                                                                                                                      & \multicolumn{4}{c}{After Mitigation}                                                                                                                                                      \\ \cmidrule(r){3-6} \cmidrule(r){7-10}
                         & & $\hat{c}_i$& $g_i$& $\lambda_i$ & Profit$^*_i$&$\hat{c}_i$ &$g_i$ &$\lambda_i$ & Profit$_i$ \\ \hline
                        \multirow{3}{*}{Non-Strategic}                            & $A$    & 20      & 23     & 20      & 0        & -       & -       & -      & 0       \\
                                                                          & $B$    & 20      & 29     & 20      & 0        & -       & -       & -      & 0       \\
                                                                          & $C$    & 20      & 23     & 20      & 0        & -       & -       & -      & 0       \\
\multirow{3}{*}{\begin{tabular}[c]{@{}c@{}}Mitigation-\\ Unaware\end{tabular}}                                                   & $A$    & 100     & 21     & 100     & 1680     & 20      & 23      & 20     & 0       \\
                                                                          & $B$    & 20      & 30     & 180     & 4800       & -       & 29       & 20      & 0       \\
                                                                          & $C$    & 20      & 24     & 20      & 0        & -       & 23       & 20      & 0       \\
\multirow{3}{*}{Conduct-Aware} & $A$    & 40      & 21     & 40      & 420      & -       & -       & -      & 420     \\
                                                                          & $B$    & 20      & 30     & 60      & 1200     & -       & -       & -      & 1200     \\
                                                                          & $C$    & 20      & 24     & 20      & 0        & -       & -       & -      & 0       \\
\multirow{3}{*}{Impact-Aware}  & $A$    & 40      & 21     & 40      & 420      & -       & -       & -      & 420     \\
                                                                          & $B$    & 20      & 30     & 60      & 1200     & -       & -       & -      & 1200\\ 
                                                                          & $C$    & 20      & 24     & 20      & 0        & -       & -       & -      & 0       \\\hline
                        \end{tabular}
        }
        \vspace{-3mm}
\end{table}

\section{6-Bus Case Study} \label{secA:6bus}
\setcounter{equation}{0}
\setcounter{figure}{0}   
\setcounter{table}{0}   

The case study with a 6-bus test system is conducted to present a more realistic bidding scenario accounting for multiple agents and multi-block offer curves. 
Let the subscript $b$ represent the generation offer block; then the submitted offer curves need to satisfy the nondecreasing constraints:
\begin{subequations}
\begin{align}
               & \hat{c}_{ib}  \ge 0,  \quad\forall i \in \Omega_{G}, b=1 \label{eq:offer1}\\
           & \hat{c}_{ib}\ge \hat{c}_{i(b-1)}  , \quad \forall i \in \Omega_{G}, \forall b \ge2. \label{eq:offer2}
\end{align}
\end{subequations}

\begin{figure}[t]
    \centering
    \includegraphics[width=0.65\linewidth]{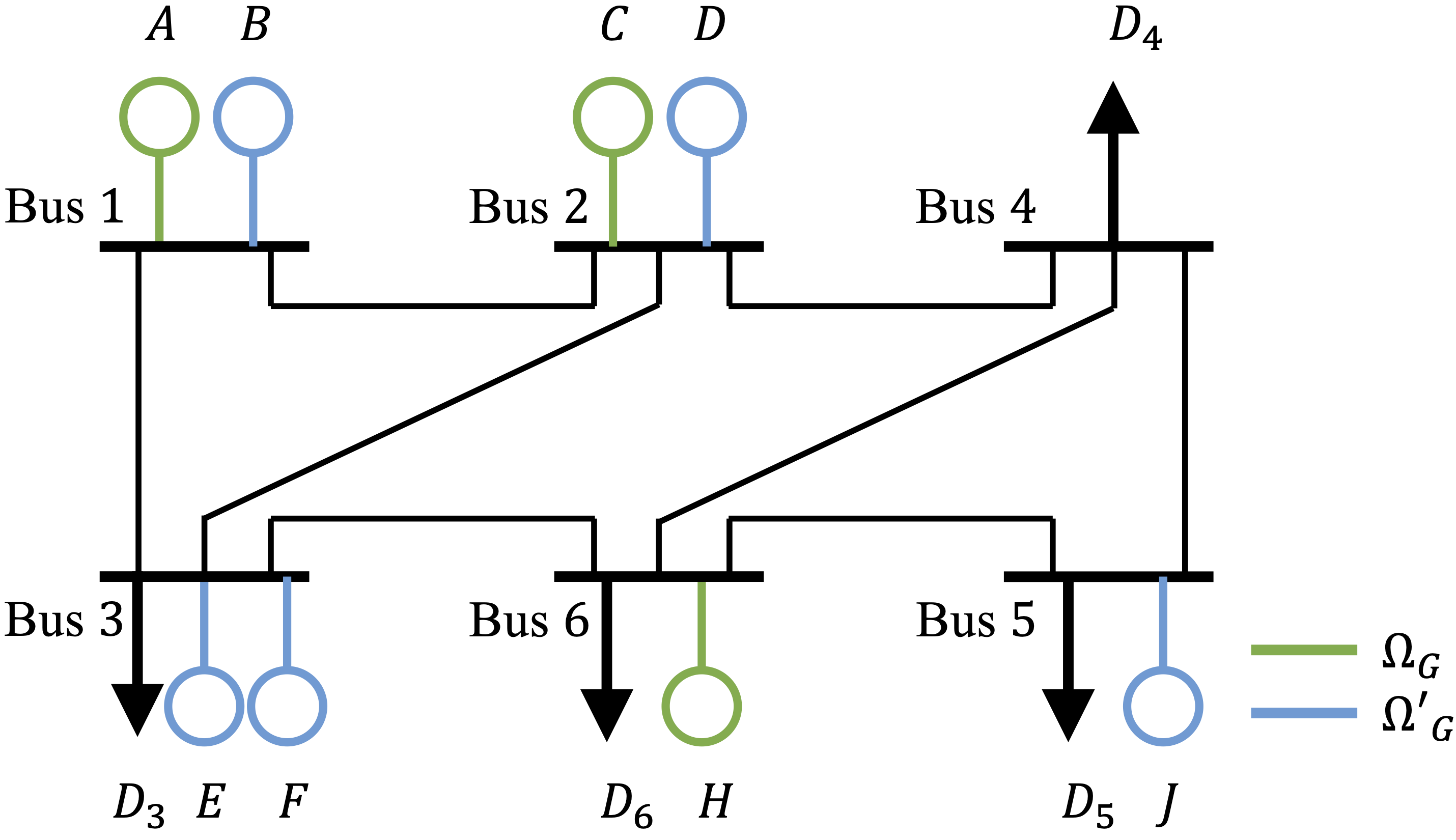}
    \vspace{-2mm}
    \caption{6-bus test system.}
    \label{fig:6bus}
    \vspace{-3mm}
\end{figure}

\begin{figure}[t]
    \centering
    \includegraphics[width=0.9\linewidth]{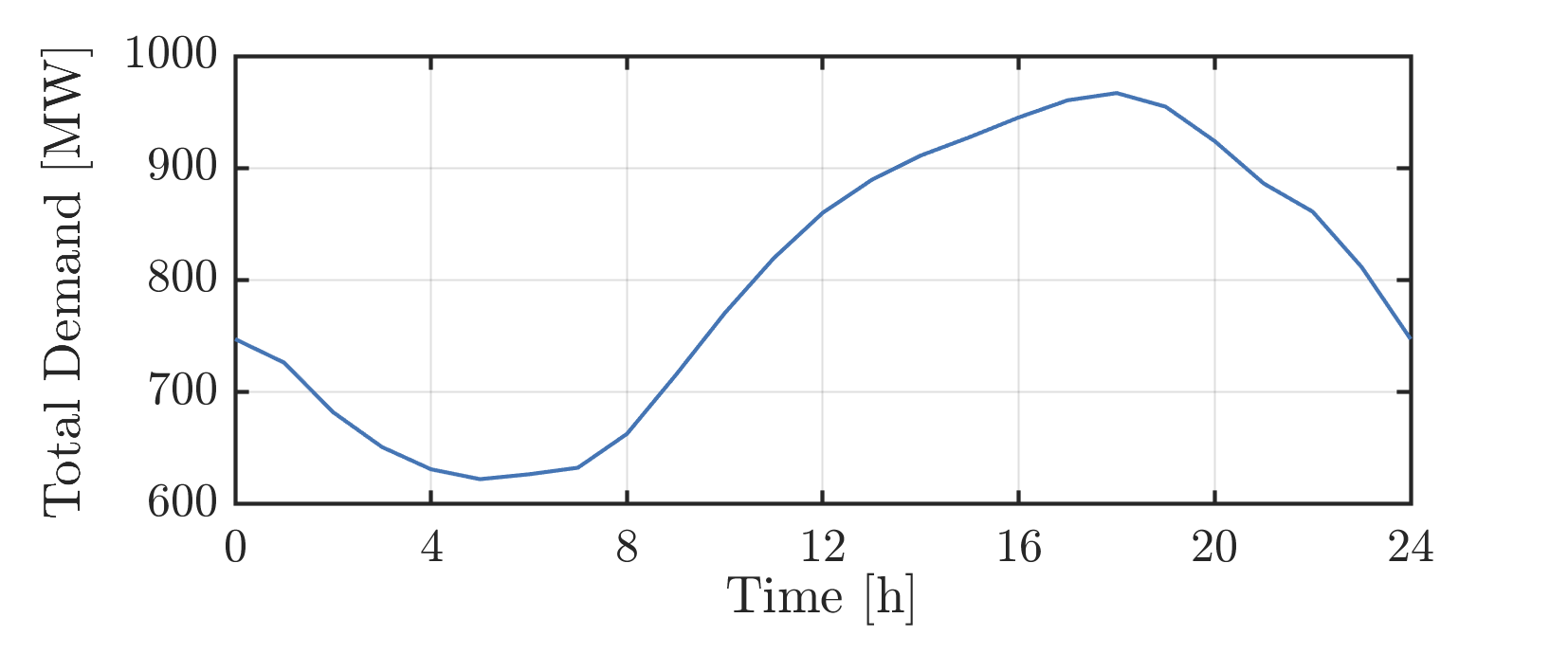}
    \vspace{-2mm}
    \caption{Demand profile.}
    \label{fig:demand}
    \vspace{-3mm}
\end{figure}

The system topology is displayed in Fig.~\ref{fig:6bus}. GenCo $G$ is considered a strategic player. The transmission line limit is 230 MW. Total demand is given in Fig.~\ref{fig:demand}. The demand share at different buses is set at $[ 0, 0, 0.19, 0.27, 0.27, 0.27]$. Hence, the test system is separated into two interconnected areas with two tie-lines: the left-hand area dominated by generation and the right-hand area dominated by load.  Table~\ref{table:6bus_set} provides a detailed configuration of the units. The offer curve of each unit is composed of four blocks.
 
\begin{table}[t]
        \caption{Data for Generating Units in the 6-Bus Test System \cite{Ruiz2009}} \label{table:6bus_set} 
        \vspace{-2mm}
        \centering 
        \setlength\tabcolsep{3pt}{
        \begin{tabular}{ccccccccc}
                \hline
                Unit                     & $A$ & $B$ & $C$ & $D$ & $E$ & $F$ & $H$ & $J$ \\ \hline
                GenCo                    & $G$ & $G'$ & $G$ & $G'$ & $G'$ & $G'$ & $G$ & $G'$ \\
                $\overline{G}_i$ {[}MW{]}    & 155   & 350   & 100   & 197   & 155   & 197   & 197  & 155    \\
                $\overline{G}_{i1}$ {[}MW{]} & 54.25 & 140.00& 25.00 & 68.95 & 54.25 & 68.95  & 68.95 & 54.25\\
                $\overline{G}_{i2}$ {[}MW{]} & 38.75 & 97.50 & 25.00 & 49.25 & 38.75 & 49.25  & 49.25 & 38.75\\
                $\overline{G}_{i3}$ {[}MW{]} & 31.00 & 52.50 & 20.00 & 39.40 & 31.00 & 39.40   & 39.40  & 31.00\\
                $\overline{G}_{i4}$ {[}MW{]} & 31.00 & 70.00 & 20.00 & 39.40 & 31.00 & 39.40   & 39.40  & 31.00\\
                $c_{i1}$ {[}\$/MWh{]}        & 9.92  & 19.20 & 18.60  & 10.08 & 9.92  & 10.08   & 10.08 & 9.92\\
                $c_{i2}$ {[}\$/MWh{]}        & 10.25 & 20.32 & 20.03 & 10.66 & 10.25 & 10.66  & 10.66 & 10.25\\
                $c_{i3}$ {[}\$/MWh{]}        & 10.68 & 21.22 & 21.67 & 11.09 & 10.68 & 11.09  & 11.09& 10.68 \\
                $c_{i4}$ {[}\$/MWh{]}        & 11.26 & 22.13 & 22.72 & 11.72 & 11.26 & 11.72  & 11.72 & 11.26\\ \hline
        \end{tabular}}
        \vspace{-3mm}
\end{table}

\begin{figure}[t]
\centering
\subfigure[Profits of GenCo $G$ ]{\includegraphics[width=0.9\linewidth]{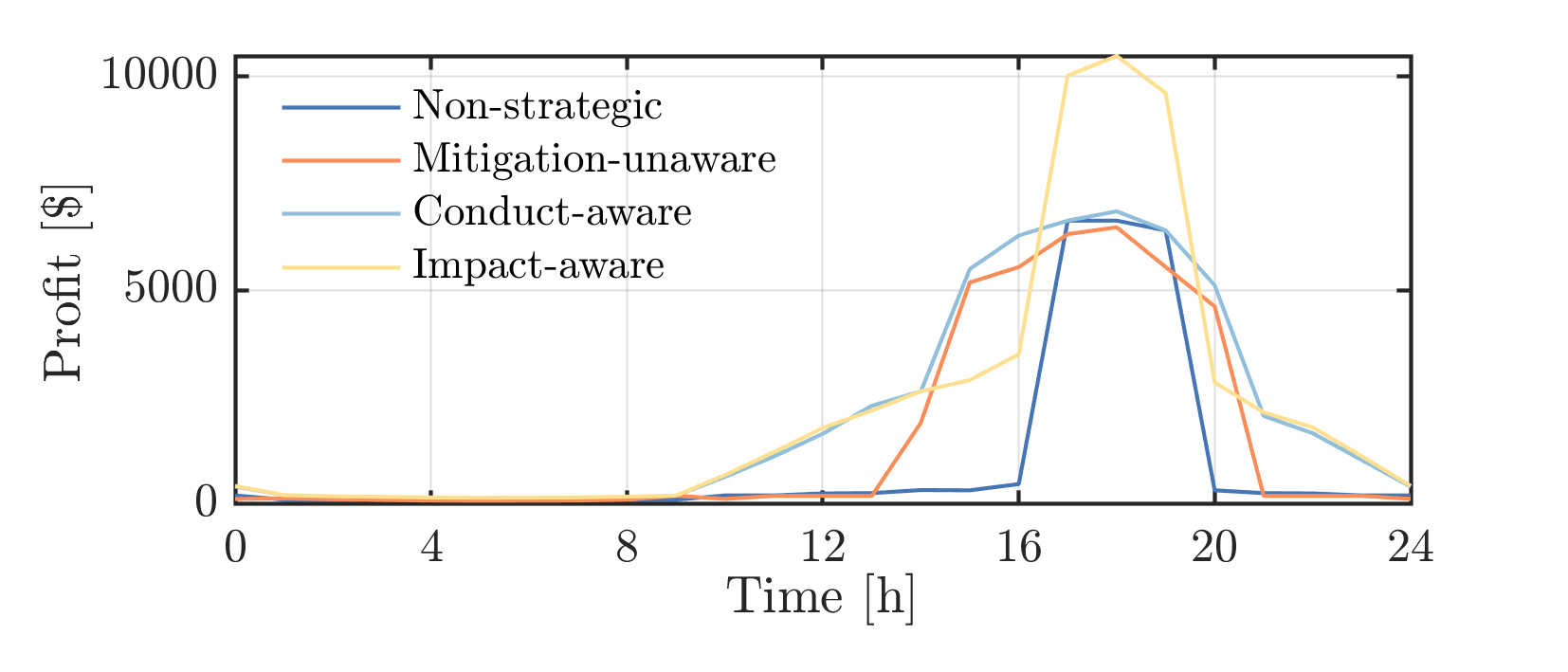} 
\label{fig:profit}} 
\subfigure[Social Welfare]{\includegraphics[width=0.9\linewidth]{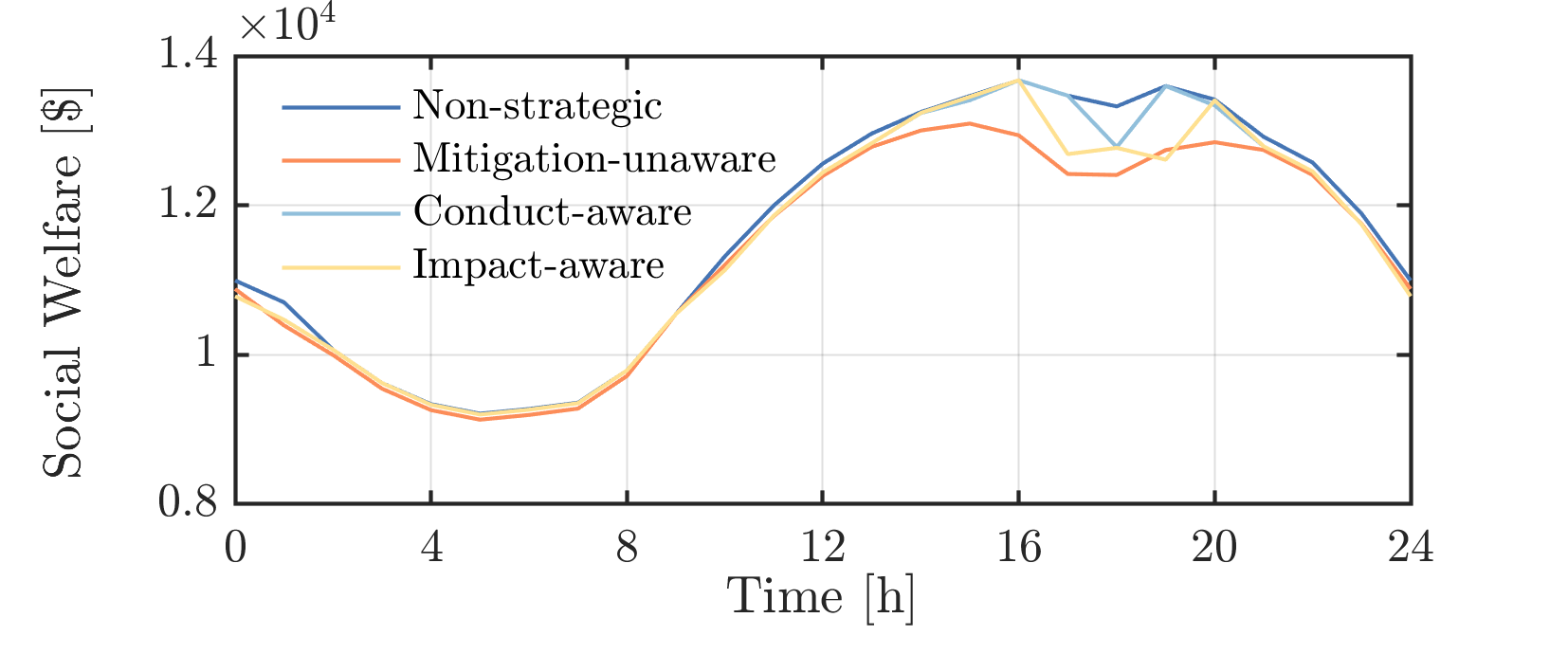} 
\label{fig:sw}}
\vspace{-2mm}
\caption{Clearing results in the congested network with the 6-bus test system.} 
\label{fig:6busresults}
\vspace{-3mm}
\end{figure}

Fig.~\ref{fig:6busresults} presents the clearing results and corresponding social welfare when GenCo $G$ uses different bidding strategies. From Fig.~\ref{fig:profit}, it can be seen that there is no significant difference among the strategies during 1:00\textendash9:00, and the reason is that the strategic GenCo lacks the driving force for exercising market power neither from the capacity limit nor the network congestion. The profits with a non-strategic bidding strategy are nonzero since the marginal cost differs at different generating levels (blocks), and units are cleared with uniform prices. 
The profit markup grows during 9:00\textendash16:00 and drops after 20:00 with the changes in demand, indicating the strong correlation between the potential of market power exercising and the demand in the system. Over these time slots, it is apparent that the conduct-aware bidding strategy outperforms non-strategic and mitigation-unaware strategies while impact-aware bidding occasionally gives a lower profit compared to the mitigation-unaware strategy, \textit{e.g.}, during 15:00\textendash16:00. It suggests the conservativeness but reliability of conduct-aware bidding and the uncertainty of impact-aware bidding associated with the estimate of clearing prices.
However, during the peak hours, \textit{i.e.}, 17:00\textendash19:00, the impact-aware bidding strategy brings high profits by catching scarcity-driven price increases in the system. Thus, the accuracy of estimating the market clearing prices is a crucial challenge for impact-aware bidding.
Another interesting aspect during this period is the occasions when the profits from mitigation-unaware bidding are lower than non-strategic bidding. This phenomenon is brought about by the aggressive bidding attempts of the strategic GenCo; once parts of the offers are mitigated, the combined profits from its units fall beyond expectations, which addresses the necessity for considering the risk of mitigation.
In terms of social welfare, as shown in Fig.~\ref{fig:sw}, it changes overall along with demand and is maximized with non-strategic bidding as a gesture of perfect competition. However, strategic bidding behaviors, either mitigation-aware or -unaware, will distort the market and render a decrease in social welfare.

Table~\ref{table:6bus} exhibits the clearing results of GenCo $G$ using the conduct-aware bidding strategy at 14:00.
Lower-priced units in the left-hand area, \textit{e.g.}, Unit $E$ and Unit $F$, first get dispatched. Line 2\textendash4 is then congested, which gives the units in the load-dominated right-hand area an ideal opportunity to exercise market power.
Considering the strategic GenCo $G$, for Unit $A$, the first two blocks are offered at the marginal cost and get fully dispatched. The offer prices for the third and fourth blocks are set to maximize the clearing price level in the left-hand area and, hence, the right-hand area. Such an offer strategy helps increase the profit of Unit $A$ and Unit $H$ as well, eventually, the GenCo.  
Unit $C$ is not cleared due to the relatively higher marginal cost.
Similarly, for Unit $H$, the last block is offered at the conduct-aware level, \textit{i.e.}, twice the true cost, to increase its profit and avoid mitigation. Note that, as a result of such strategic behavior, the clearing price for Unit $G$ is \$29.15/MWh, which is higher than that for Unit $H$. 
Given the multi-block bidding policy, a strategic GenCo can design subtle offer strategies to maximize its profit.

\begin{table}[ht]
        \caption{Clearing Results of GenCo $G$ Using Conduct-Aware Bidding at 14:00 } \label{table:6bus} 
        \vspace{-2mm}
        \centering 
        \setlength\tabcolsep{5pt}
        {
        \begin{tabular}{cccccccc}
        \hline
Unit  & $\hat{c}_{i1}$ & $\hat{c}_{i2}$ & $\hat{c}_{i3}$ & $\hat{c}_{i4}$ & $g_i$  & $\lambda_i$ & Profit$_i$ \\\hline
$A$ & 9.92           & 10.25          & 15.20          & 15.20          & 119.08 & 15.20       & 595.94     \\
$C$ & 18.60          & 20.03          & 21.67          & 22.72          & 0.00   & 11.09       & 0.00       \\
$H$ & 10.08          & 10.66          & 11.09          & 23.44          & 157.60  & 23.44       & 2037.16    \\\hline
\end{tabular}}
        \vspace{-3mm}
\end{table}




\end{document}